\documentclass[12pt,reqno]{amsart} 

\usepackage{graphicx,amssymb,colordvi,textcomp,latexsym,mathtools,enumitem}
\usepackage[hidelinks]{hyperref}
\usepackage{multirow}
\usepackage{tabularx}
\usepackage{booktabs}
\usepackage{float}
\usepackage{prettyref}
\usepackage{chngcntr}
\counterwithout{equation}{subsection}
\counterwithout{equation}{section}

\def\ddefloop#1{\ifx\ddefloop#1\else\ddef{#1}\expandafter\ddefloop\fi}
\def\ddef#1{\expandafter\def\csname 
	bb#1\endcsname{\ensuremath{\mathbb{#1}}}}
\ddefloop ABCDEFGHIJKLMNOPQRSTUVWXYZ\ddefloop
\def\ddef#1{\expandafter\def\csname 
	#1\endcsname{\ensuremath{\mathbf{#1}}}}
\ddefloop abcdefghijklmnopqrstuvwxyz\ddefloop
\def\ddefloop#1{\ifx\ddefloop#1\else\ddef{#1}\expandafter\ddefloop\fi}
\def\ddef#1{\expandafter\def\csname 
	b#1\endcsname{\ensuremath{\mathbb{#1}}}}
\ddefloop ABCDEFGHIJKLMNOPQRSTUVWXYZ\ddefloop
\def\ddef#1{\expandafter\def\csname 
	c#1\endcsname{\ensuremath{\mathcal{#1}}}}
\ddefloop ABCDEFGHIJKLMNOPQRSTUVWXYZ\ddefloop
\def\ddef#1{\expandafter\def\csname 
	f#1\endcsname{\ensuremath{\mathfrak{#1}}}}
\ddefloop ABCDEFGHIJKLMNOPQRSTUVWXYZabcdefghjklmnopqrstuvwxyz\ddefloop
\def\ddef#1{\expandafter\def\csname 
	h#1\endcsname{\ensuremath{\widehat{#1}}}}
\ddefloop ABCDEFGHIJKLMNOPQRSTUVWXYZ\ddefloop
\def\ddef#1{\expandafter\def\csname 
	hc#1\endcsname{\ensuremath{\widehat{\mathcal{#1}}}}}
\ddefloop ABCDEFGHIJKLMNOPQRSTUVWXYZ\ddefloop

\newcommand{\balpha}{\boldsymbol{\alpha}}

\newcommand{\bbeta}{\boldsymbol{\beta}}

\newcommand{\tens}{\otimes{}}
\newcommand{\defeq}{\coloneqq}
\DeclarePairedDelimiter{\brk}{[}{]}

\DeclarePairedDelimiter{\prn}{(}{)}
\DeclarePairedDelimiter{\nrm}{\|}{\|}
\DeclarePairedDelimiter{\inner}{\langle}{\rangle}

\newcommand{\symp}{\mathrm{Sym}}
\newcommand{\sympn}{\symp^n(\reals^d)}

\newcommand{\reals}{\mathbb{R}}
\newcommand{\isov}{\Gamma}

\theoremstyle{definition}  %Sets style of subsequent newtheorems to 'definition'

\newtheorem{cor}{Corollary}
\newtheorem{prop}{Proposition}
\newtheorem{fact}{Fact}

\theoremstyle{plain}
\newtheorem{rem}{Remark}
\newtheorem{theorem}{Theorem}
\newtheorem{definition}{Definition}

\usepackage{prettyref}
\newcommand{\pref}[1]{\prettyref{#1}}
\newcommand{\pfref}[1]{Proof of \prettyref{#1}}
\newcommand{\savehyperref}[2]{\texorpdfstring{\hyperref[#1]{#2}}{#2}}
\newrefformat{eq}{\savehyperref{#1}{\textup{(\ref*{#1})}}}
\newrefformat{eqn}{\savehyperref{#1}{Equation~\ref*{#1}}}
\newrefformat{lem}{\savehyperref{#1}{Lemma~\ref*{#1}}}
\newrefformat{obs}{\savehyperref{#1}{Observation~\ref*{#1}}}
\newrefformat{def}{\savehyperref{#1}{Definition~\ref*{#1}}}
\newrefformat{thm}{\savehyperref{#1}{Theorem~\ref*{#1}}}
\newrefformat{corr}{\savehyperref{#1}{Corollary~\ref*{#1}}}
\newrefformat{sec}{\savehyperref{#1}{Section~\ref*{#1}}}
\newrefformat{app}{\savehyperref{#1}{Appendix~\ref*{#1}}}
\newrefformat{ass}{\savehyperref{#1}{Assumption~\ref*{#1}}}
\newrefformat{ex}{\savehyperref{#1}{Example~\ref*{#1}}}
\newrefformat{fig}{\savehyperref{#1}{Figure~\ref*{#1}}}
\newrefformat{alg}{\savehyperref{#1}{Algorithm~\ref*{#1}}}
\newrefformat{rem}{\savehyperref{#1}{Remark~\ref*{#1}}}
\newrefformat{conj}{\savehyperref{#1}{Conjecture~\ref*{#1}}}
\newrefformat{prop}{\savehyperref{#1}{Proposition~\ref*{#1}}}
\newrefformat{table}{\savehyperref{#1}{Table~\ref*{#1}}}
\newrefformat{prob}{\savehyperref{#1}{Problem~\ref*{#1}}}
\newrefformat{claim}{\savehyperref{#1}{Claim~\ref*{#1}}}
\newrefformat{fact}{\savehyperref{#1}{Fact~\ref*{#1}}}
\newrefformat{sec}{\savehyperref{#1}{Section~\ref*{#1}}}
\newrefformat{exam}{\savehyperref{#1}{Example~\ref*{#1}}}
\newrefformat{exams}{\savehyperref{#1}{Examples~\ref*{#1}}}
\newrefformat{prob}{\savehyperref{#1}{(\ref*{#1})}}

\usepackage{setspace}

\newcommand{\RR}{\mathbb{R}}
\newcommand{\EE}{\mathbb{E}}
\newcommand{\QQ}{\mathbb{Q}}

\newcommand{\ZZ}{\mathbb{Z}}
\newcommand{\NN}{\mathbb{N}}

\newcommand{\ploss}{{\cL}} 
\renewcommand{\ker}{{\kappa}}

\newcommand{\bxi}{{\boldsymbol{\xi}}}

\newcommand{\grad}[1]{{\nabla {#1}}}

\title[Critical points for Symmetric Tensor Decomposition]{Symmetry  \& 
Critical Points for Symmetric Tensor Decomposition Problems}

\author{Yossi Arjevani and Gal Vinograd}
\address{Yossi Arjevani, The Hebrew University, Jerusalem}
\email{yossi.arjevani@gmail.com}
\address{Gal Vinograd, The Hebrew University, Jerusalem}
\email{gal.vinograd@gmail.com}

\begin{document}

\begin{abstract}

We consider the nonconvex optimization problem associated with the decomposition of a real symmetric tensor into a sum of rank-one terms. Use is made of the rich symmetry structure to construct infinite families of critical points represented by Puiseux series in the problem dimension, and so obtain precise analytic estimates on the objective function value and the Hessian spectrum. The results enable an analytic characterization of various obstructions to local optimization methods, revealing, in particular, a complex array of saddles and minima that differ in their symmetry, structure, and analytic properties. A~notable phenomenon, observed  for all critical points considered, concerns the index of the Hessian  increasing with the objective function value. 
\end{abstract}

\maketitle

\section{Introduction}\label{sec:intro}
The need to decompose a real symmetric tensor into a sum of rank-one terms arises naturally in various scientific domains \cite{comon2008symmetric,comon2006blind,de2000multilinear,smilde2005multi,mccullagh2018tensor}. Standard approaches
\cite{carroll1970analysis,harshman70,paatero1997least,tomasiuse,acar2011scalable,kolda2015numerical} involve tackling a nonconvex optimization problem associated with an order-$n$ symmetric tensor $T$ on $\RR^d$, given by  
\begin{align} \label{prob:opt}
\min_{
	\stackrel{W\in M(k,d),}{\balpha\in\reals^k}}\ploss(W,\balpha) 
	\defeq  
\nrm*{\sum_{i=1}^k \alpha_i \w_i^{\tens n} - T}^2
\end{align}
for fixed $k\in \NN$, where $M(k,d)$ denotes the set of $k\times d$-real matrices, $\w_i$ the $i$th row of a weight matrix $W\in M(k,d)$, and $\nrm{\cdot}$ a tensor norm. It has been observed 
\cite{acar2011scalable,kolda2015numerical,kolda2009tensor,tomasi2006comparison} that these approaches often fail to find valid decompositions in various parameter regimes encountered in practice. This article focuses on the geometric obstructions that arise in the associated optimization landscape, particularly saddles and minima, which contribute to such failures. Rather than striving for general results, we focus on symmetric tensors that are invariant to certain actions of subgroups of $S_m$, the symmetric group over $m$ symbols. This choice yields analytically tractable problems whose geometry exposes key obstructions to symmetric tensor decomposition—the aim of this~work.

\section{Main results} \label{sec:main_results}
After a quick review of the setting studied and key related concepts, 
the introduction continues with a detailed description of the main 
results. We conclude with a brief survey of related work and outline of the structure of the paper. Some familiarity with real algebraic geometry and group theory is assumed, but otherwise, key ideas and concepts are introduced as needed, with the exception of several \emph{irreducible representations} of $S_d$ that are explicitly referred to below (see \pref{sec:sb}) and may be safely skipped without loss of continuity. Proofs are deferred to the appendix.

We consider a natural action of $S_k\times S_d$ on the parameter space $M(k,d)$ (orthogonal): the first factor permutes rows, the second columns. Given a weight matrix $W\in M(k,d)$, the largest subgroup of $S_k\times S_d$ fixing $W$ is called the \emph{isotropy group} and is used as a means of measuring the symmetry of $W$ (see \pref{sec:inv}). For example, for $k=d$, the isotropy group  of the identity matrix is the diagonal subgroup $\Delta S_d \defeq \{(\pi,\pi)~|~\pi\in S_d\}\subseteq S_d\times S_d$. In \cite{arjevani2021symmetry}, it was shown that, if $T$ possesses certain invariance properties, $\ploss$ is $S_k\times S_d$-invariant. An intricate mechanism, connecting the symmetry of adjacent critical points of invariant functions \cite{arjevani2024symmetry}, then dictates the emergence of critical points of $\ploss$ that exhibit \emph{proper} isotropy subgroups \emph{breaking the symmetry} of $S_k \times S_d$, as is  verified through numerical and algebraic methods. In this work, this phenomenon of symmetry breaking is employed to derive sharp analytic results for the optimization landscape associated with $\ploss$. More specifically,

\begin{itemize}[leftmargin=*]
\item The invariance properties of a broad class of optimization problems are formalized, encompassing those associated with symmetric tensor decomposition and two-layer ReLU networks. These properties enable a detailed study of symmetric critical points, particularly ones exhibiting large isotropy groups.

\item We construct infinite families of critical points represented by Puiseux series in $d$ and derive sharp analytic estimates of the Hessian spectrum. The tensor norms considered are the Frobenius norm and the cubic-Gauss norm defined in (\ref{gauss_nrm}).

\item The complex array of saddles and minima so revealed exhibits various phenomena affecting the behavior of local optimization methods, as discussed in length in \pref{sec:ss}. For example, saddles with a higher objective value are seen to exhibit a larger number of descent~directions.

\end{itemize}

A general introduction to the symmetry-breaking framework is provided in \pref{sec:sb}. For reasons of exposition, the applications in this article are presented primarily for third-order symmetric tensors under the assumption $k=d$,  following for example \cite{tomasi2006comparison}. Moreover, we shall focus on symmetric tensors admitting an orthonormal decomposition, that is, $T = \sum_{i=1}^d \v_i^{\tens 3}$ where $\v_1,\dots,\v_d$ are orthonormal. The methods, however, are quite general. Orthonormal decomposability of $T$ (see related \cite{robeva2016orthogonal}) is not required; other choices of order, dimensionality, rank and target tensors with possibly lesser symmetry may be considered.

\subsection{Global minimizers} \label{sec:global}
We begin by providing a simple characterization of the global minimizers of 
$\ploss$. Our study of families of critical points is essentially invariant to orthogonal transformations of the domain of $\ploss$, and so throughout we assume (without loss of generality, \pref{sec:generalities}) that the target tensor is $T_\e \defeq \sum_{i=1}^d \e_i^{\tens 3}$, where $\e_i$ denotes the $i$th unit vector. Moreover, since $n$ is odd, all $\alpha_i$ may be absorbed into $\w_i$ and are thus set to 1. Given a weight matrix $W\in M(d,d)$, the action of $S_d\times S_d$, perturbing rows and columns of $W$, typically yields additional points 
constituting $(S_d\times S_d)W$, the \emph{group orbit} of~$W$. For 
our choice $T=T_\e$, the objective loss $\ploss$ is $S_d\times 
S_d$-invariant (see \pref{sec:inv}) and so points on the orbit of a 
critical point are also critical and share the same loss and Hessian 
spectrum, see generally \cite{golubitsky2012singularities}. Since 
$\ploss\ge 0$ and $\ploss(I_d) = 0$, the identity matrix $I_d$, as well 
as all points lying on $(S_d\times S_d) I_d$, are global minimizers. That 
the converse holds is a consequence of Kruskal's criterion for 
\emph{identifiability}~\cite{kruskal1977three}.
\begin{prop}\label{prop:global}
Independently of the choice of the tensor norm, 
a point $W\in M(d,d)$ is a global minimizer of $\ploss$ iff~$W\in 
(S_d\times S_d) I_d$, the set of all $d\times d$ permutation matrices. 
In particular, the set of global minimizers and the set of points defining 
valid decompositions for $T_\e$ are~identical.
\end{prop}
Versions of the proposition involving the action of the \emph{Hyperoctahedral} group hold for tensors of order $n\ge4$ and may be proven by an extension to Kruskal's criterion given in \cite{sidiropoulos2000uniqueness} (see \cite[Prop. 4.14]{ArjevaniField2020} 
for an analog result for ReLU).

By \pref{prop:global}, issues concerning  the existence of global solutions for (\ref{prob:opt}) do not arise \cite{de2008tensor}; a global solution always exists. Moreover, the number of global minimizers is equal to the \emph{length }of the orbit $|(S_d \times S_d)I_d| = d!$. The large number of global minimizers accounts for, in part, the complex array of critical points occurring for $\ploss$. In the general case, the orbit-stabilizer theorem states that the smaller the isotropy group, the larger the group orbit. Below we indicate orbit lengths at appropriate points. 

Although the set of global minimizers is discrete, this need not be the case for other families of critical points. For example, a continuous of family of critical points, studied in some depth below, is given by 
$\fC_{5,t}(d) \defeq I_{d-2} \oplus \begin{bsmallmatrix} (1-t^3)^{1/3} & 
0\\ t & 0\end{bsmallmatrix}$. 
\begin{prop}\label{prop:proper}
For all $t\in\RR$ and $d\ge2$, $\fC_{5,t}(d)$ is a critical point having $\ploss(\fC_{5,t}(d))=1$. In particular, the level set $\ploss^{-1}(\{1\})$ is not bounded and so $\ploss$ is not a proper map.
\end{prop}
The family of critical points $\fC_{5,0}$ in fact defines saddles. The proof is based on the curve selection (CSL) lemma, used in the framework of tangency arcs; see \cite{arjevani2024symmetry}.

\subsection{Symmetric Saddles} \label{sec:ss}
Contrarily to global minimizers, the set of saddle points 
exhibits a highly complex and intricate structure. Below, we 
consider two tensor norms widely used in the 
literature (see \pref{sec:inv}): the Frobenius norm $\|\cdot\|_F$ 
and the \emph{cubic-Gaussian} norm $\|\cdot\|_\cN$ 
defined by
\begin{align}\label{gauss_nrm}
	\|S-T\|_{\cN} \defeq \EE_{\x\sim  \cN(0,I_d)}\brk{  \inner{S-T, 
			\x^{\tens 3}}^2
	},
\end{align}
with $\cN(0,I_d)$ denoting the standard multivariate Gaussian distribution, see \cite{arjevani2021symmetry,render2008reproducing} and references therein. The objective function will be respectively denoted by $\ploss_F$ and $\ploss_\cN$. A consequence of the study of the two different norms is that phenomena which may otherwise be deemed norm-specific can be examined and verified in separate settings, see for example (\ref{num:conj}) below.  

\subsubsection{The Frobenius norm}
We begin by considering families of critical points $\{W(d)\in M(d,d)\}_{d\in \NN}$ occurring for the Frobenius norm. The entries of each family are expressed in terms of Puiseux series in $d$, giving a critical point of $\ploss_F:M(d,d)\to\RR$ for every integer $d\in \NN$.

\begin{theorem} \label{thm:for}
The sequences of weight matrices $\fC_i(d)\in M(d,d),~i\in \{I,0,1,2,3,4,5\}$,
given in \pref{sec:sb}, define families of critical points of $\ploss_F$ for 
the target tensor $T_\e$. 
The loss, the Hessian spectrum and the orbit length are 
given in 
\pref{table:sdTsdevs} and 
\pref{table:devs}.
\end{theorem}

\begin{table}[ht]
	\begin{center}
		\begin{tabular}{|l|l||c|c|c|c| }\cline{3-6}
\multicolumn{2}{l||}{}
&		
$\fC_0$&$\fC_1$ & \hspace*{0.05in}$\fC_2$& $\fC_I$\\\hline\hline
\multicolumn{2}{|l||}{Loss}&		
$d$&$d-1/d$&  $d-1/d$&$0$\\\hline
\multicolumn{2}{|l||}{Isotropy}&		
$S_d\times S_d$ &$S_d\times S_d$ & $\Delta S_d$ & $\Delta S_d$\\\hline
\multicolumn{2}{|l||}{Orbit length} &   $1$& $1$&$d!$&$d!$\\\hline
	Comp.  &Multiplicity 
	\hspace*{0.05in}\\\hline
	$\ft$& $1$ & $0$&$18/d$  & $18/d$&$18$\\
	&&&& $12/d$&$6$\\
	\hline
	$\fs_d$ & $d-1$ & $0$&$0$& $12/d$ &$18$ \\
&&$0$&$-6/d$&$0$&$6$\\
	&&&& $-6/d$&$6$\\
	\hline
	$\mathfrak{s}_{d}\tens \mathfrak{s}_{d}$ & $(d-1)^2$ & 
	$0$&$-12/d$  &&\\
	\hline
	$\fx_d$ & $(d-1)(d-2)/2$ & &&	$-12/d$ &$6$ \\
	\hline
	$\fy_d$ & $d(d-3)/2$ & &	&$-12/d$&$6$  \\
	\hline
	
\end{tabular}
\end{center}
\caption{
The loss and the Hessian spectrum given in $o(d^{-1})$-terms for the 
families of critical points $\fC_1$ and $\fC_2$, and in precise terms 
for $\fC_I$ giving the identity matrices (global minimizers, see 
\pref{prop:global}) and for $\fC_0$ the zero matrices. The families $\fC_1$ and 
$\fC_2$ define saddles for sufficiently large $d$. The family $\fC_0$ (as well as $\fC_5$ in \pref{table:devs} below) also defines saddles, but due to the singularity of the Hessian, the proof requires going beyond second-order analysis; see \cite{arjevani2024symmetry}. The orbit length indicates the number of different families of critical points that occur by the action of $S_d\times S_d$, having in particular identical loss and Hessian spectrum which is listed by the respective irreducible representation of  $S_d$, \pref{sec:sb}. 
}
\label{table:sdTsdevs}
\end{table}

\begin{table}
\makebox[0pt]{\begin{tabular}{|l||l|l|l|l|l|l|} \hline
		&
		\multicolumn{2}{|c|}{$\fC_3$}&\multicolumn{2}{|c|}{$\fC_4$} 
		&\multicolumn{2}{|c|}{ $\fC_5\defeq \fC_{5,0}$}\\\hline\hline
		Isotropy &\multicolumn{2}{|c|}{
		$\Delta (S_{d-1}\times S_1)$} &
		\multicolumn{2}{|c|}{$\Delta(S_{d-2}\times S_2)$} & 
		\multicolumn{2}{|c|}{$\Delta 
		(S_{d-2}\times 			S_1^2)$}\\\hline
		Loss&\multicolumn{2}{|c|}{$\quad\quad  
		d-1$}&\multicolumn{2}{|c|}{$3/2$}&  
		\multicolumn{2}{|c|}{$1$}\\\hline
		Orbit length&\multicolumn{2}{|c|}{$\quad  
d\cdot d!$}&\multicolumn{2}{|c|}{$d(d-1)d!/2$}&  
\multicolumn{2}{|c|}{$d(d-1)d!$}\\\hline\hline
&Eigenvalue & $\#$ &Eigenvalue & 
$\#$&Eigenvalue & 
$\#$\\\hline
&$18$&$1$&$18$&$d-2$&$18$&$d-1$\\
&$6$&$d-1$&$9$&$1$&$6$&$(d-1)^2$\\
&$18/d$&1&$6$&$d^2-3d+2$&$0$&$d$\\
&$12/d$&$d-1$&$3$&$d-2$&&\\
&$6/d$&1&$0$&$d-1$&&\\
&$0$&$2d-4$&$-3$&$1$&&\\
&$-6/d$&$d-2$&$-6$&$1$&&\\
&$-12/d$&$d^2-5d+5$&&&&\\\hline
	\end{tabular}}
\caption{
The loss and the Hessian spectrum given in $o(d^{-1})$-terms for the family of critical points $\fC_3$ and in precise terms for the families of critical points $\fC_4$ and~$\fC_5$. Similarly to $\fC_0$ appearing in \pref{table:sdTsdevs}, it is not possible to determine whether $\fC_5$ defines a family of saddles or minima by inspecting the Hessian spectrum alone as  the Hessian is positive semi-definite and singular. In~\cite{arjevani2024symmetry}, the existence of real analytic descent curves for $\fC_5$ is proved using the CSL~\cite{milnor61singular}, showing in particular that $\fC_5$ defines saddles.
}
\label{table:devs}
\end{table}

In \cite{cai2022gradient}, several families of minima of $\ploss_F$ 
are described. However, only the $k=1$ case is considered, and the Hessian spectrum is provided only for certain local minima and with limited accuracy. In contrast, the symmetry breaking approach applies in principle to any choice of $k$ and $d$, and yields an analytic description of the Hessian spectrum (as well as the critical value) at both local minima and saddle points. We note that although the Hessian has $d^2$ eigenvalues, the number of \emph{distinct }eigenvalues is fixed and independent of the dimension; a consequence of the symmetry structure, namely, the multiplicity of irreducible representations occurring in the respective  \emph{isotypic decomposition} being independent of the dimension (\pref{sec:sb}). 

The saddle points given in \pref{thm:for} behave rather differently in terms of the dependence of their objective value and Hessian spectrum on the dimension, cf. \cite{ArjevaniField2021}. While the Hessian spectrum at $\fC_0, \fC_1$ and $\fC_2$ converges to zero, some of the positive eigenvalues of $\fC_3$ and $\fC_4$ are bounded away from zero. In fact, for $\fC_4$, except for two negative eigenvalues, $\Theta(d^2)$ of the eigenvalues are strictly positive, and so despite being saddles, critical points belonging to  $\fC_4$ may be hard to escape via gradient-based methods. In contrast, at $\fC_3$, $\Theta(d^2)$ of the Hessian eigenvalues are negative (of the order of $1/d$). Comparing the objective value, $\ploss_\cF$ evaluated at either $\fC_0, \fC_1, \fC_2$ or $\fC_3$ scales as $\Theta(d)$. However, critical values given by saddles need not grow linearly with~$d$. For example, at $\fC_4$, the critical value is fixed and equal to $3/2$. 

Our findings suggest a direct connection between the function value at non-zero critical points and the \emph{index}, that is, the number of negative Hessian eigenvalues, at a (not necessarily \emph{hyperbolic}) critical point: the larger the critical value, the larger the index---a desirable property for successful use of local optimization methods. Quantitatively, for a general family critical points $\fC(d)$, this may be expressed in the form 
\begin{align}\label{num:conj}
\frac{\ploss(\fC(d))}{d} \approx 
\frac{\mathrm{Index}(\fC(d))}{d^2},
\end{align}
or, rather, in a more general form replacing second-order descent 
directions, counted by $\mathrm{Index}(\cdot)$, with ones corresponding to second- \emph{or} higher-order derivatives. The above phenomenon (indeed, its quantitative version  (\ref{num:conj})) also occurs for the cubic-Gaussian norm (\pref{table:normal_sdTsdevs}  below) and for two-layer ReLU networks 
\cite{ArjevaniField2021} (cf. \cite[Remark 2]{ArjevaniField2022annihilation} on the {expected }initial value growing linearly with $d$ under Xavier initialization). It is therefore natural to conjecture that the above, or a {qualitatively} similar, relation between the function value at critical points and the index applies more broadly to non-convex functions exhibiting symmetry breaking phenomena of the nature considered in this work.

\subsubsection{{The cubic-Gaussian norm}}
Although the geometry of the associated optimization landscape depends on the choice of the tensor norm, by \pref{prop:global}, the set of global minimizers giving valid tensor decomposition is identical for all norms. The critical set, however, differs quite significantly. We illustrate this for the cubic-Gaussian norm (\ref{gauss_nrm}).

\begin{theorem} \label{thm:gauss}
The sequences of weight matrices $\fD_I(d), \fD_1$ and $\fD_2$, described in \pref{sec:sb}, define families of critical points of $\ploss_\cN$ for the target tensor $T_\e$. The loss, the Hessian spectrum and the orbit length are given in  \pref{table:normal_sdTsdevs}.
\end{theorem}

The loss at $\fD_I,\fD_1$ and $\fD_2$ scale similarly to ones seen for the Frobenius norm. However, the Hessian spectrum behaves rather differently. For example, the Hessian spectrum at the families of minima $\fD_I$ is extremely skewed having $\Theta(d^2)$ positive eigenvalues concentrated near zero and $\Theta(d)$ eigenvalues growing linearly with $d$. In addition, modulo $\Theta(d)$ eigenvalues, the spectrum of the families of third-order {saddles} $\fD_2$  and $\fD_I$ essentially agree. The two phenomena: (a) highly-skewed Hessian spectra (b) largely agreeing at certain critical points, have also been observed for two-layer ReLU networks \cite{arjevanifield2020hessian,ArjevaniField2021}. Lastly,  although $\fC_I$ and $\fD_I$ represent the same family of global minima, the Hessian spectra differ significantly: the former is fixed, the latter depends on $d$. This further emphasizes  the dependence of the geometry of the associated optimization landscape on the choice of the tensor~norm.

\begin{table}[h]
	\begin{center}
		\begin{tabular}{ |c|c||c|c|c|c| } \cline{3-5}
			\multicolumn{2}{c||}{} &
			\hspace*{0.05in}$\fD_1$ & 
			$\fD_{I}$&	$\fD_{2}$\\\hline\hline
			\multicolumn{2}{|c||}{Isotropy}&$S_d\times S_d$&$\Delta 
			S_d$&$\Delta 
			(S_{d-1}\times S_1)$\\\hline
			\multicolumn{2}{|c||}{Loss} &   $\frac{48 d}{5} - 
			\frac{36}{5} - 
			\frac{12}{5 d}$& 0&$6+18/d$\\\hline
			\multicolumn{2}{|c||}{Multiplicity} &   $1$& $d!$&$d\cdot 
			d!$\\\hline
			Comp.& Degree&\multicolumn{3}{c}{}
			%			Loss&		
			\\\hline
			$\ft$& 1 & $\frac{162 \cdot {75}^{\frac{1}{3}} 
			}{5}{d}^{\frac{1}{3}}  + \frac{144\cdot {75}^{\frac{1}{3}}}{5 
			}d^{\frac{-2}{3}}$  & 108&0\\
			&&&$18d+180$&0\\
			&&&&108\\
			&&&&$18d+18$\\
			&&&&$18d+162$\\
			\hline
			$\mathfrak{s}_{d}$& $d-1$ & 0& 108&\\
			&&$\frac{18 \cdot {75}^{\frac{1}{3}} }{5}{d}^{\frac{1}{3}} - 
			\frac{32 \cdot {75}^{\frac{1}{3}}}{5 
			}d^{\frac{-2}{3}}$&36&\\
			&&& $18d+180$&\\
			\hline
			$\mathfrak{s}_{d-1}$ & $d-2$ && &0\\
			&&&&36\\
			&&&&36\\
			&&&&108\\
			&&&&$18d+162$\\
			\hline
			$\mathfrak{s}_{d}\tens \mathfrak{s}_{d}$ & $(d-1)^2$ & 
			$- \frac{72 \cdot {75}^{\frac{1}{3} }}{25}{d}^{\frac{1}{3}} 
			- 
			\frac{304 \cdot {75}^{\frac{1}{3}}}{25 }d^{\frac{-2}{3}}$  
			&&\\
			\hline
			$\fx_d$ & $(d-1)(d-2)/2$ &  &36&\\
			\hline
			$\fx_{d-1}$ & $(d-2)(d-3)/2$ &  &&36\\
			\hline
			$\fy_d$ & $d(d-3)/2$ & &36&\\
			\hline
			$\fy_{d-1}$ & $(d-1)(d-4)/2$ & &&36\\
			\hline
		\end{tabular}
	\end{center}
	\caption{\label{table:normal_sdTsdevs}
The loss and the Hessian spectrum of three families of 
critical points given in \pref{sec:rank_for}, with $\fD_I$ defining global solutions and $\fD_1$ saddles. The Hessian spectrum is given modulo $o(d^{-\frac{2}{3}})$-terms and is listed by the respective irreducible representation of  $S_d$ (\pref{sec:sb}).}
\end{table}

\subsection{Related work}
We relate our results to the existing literature.

\subsubsection{Symmetric tensor decomposition.}
The tensor decomposition problem considered in this work, also
referred to as the polyadic decomposition, was introduced 
under the names CANDECOMP and PARAFAC by respectively
\cite{carroll1970analysis} and \cite{harshman1970foundations}.  
Standard approaches involving nonconvex optimization formulation include
\cite{carroll1970analysis,
	harshman70,paatero1997least,tomasiuse,
	acar2011scalable,kolda2015numerical}.
See \cite{kolda2009tensor} for a survey of methods and 
applications. As noted in the introductory comments, it has been observed \cite{acar2011scalable,kolda2015numerical,kolda2009tensor,tomasi2006comparison} that these approaches have a relatively high likelihood of failure across various parameter regimes encountered in practice. Theoretical guarantees are available in some cases \cite{cai2022gradient,liu2022symmetric,jiang2015tensor}, but none of the latter  directly addresses the optimization landscape for $k>1$. A  computationally intensive yet provable method, assuming $T$ has rank $O(n^{\frac{d-1}{2}})$, was presented in \cite{nie2017generating}. For $n=3$ and $k \leq d$, a theoretically convenient method was derived based on simultaneous diagonalization of matrix slices \cite{harshman70}.   For $k=4$ and $k = O(n^2)$, the work \cite{de2007fourth} introduced a provable algorithm using matrix eigendecompositions, which was later robustified using ideas from the sums-of-squares hierarchy in \cite{hopkins2019robust}. In \cite{kileel2019subspace}, a tensor power method, constructed from a matrix flattening of $T$, was applied to sequentially extract the components 
$\w_i$.

\subsection{Organization of the paper.} In \pref{sec:preliminaries}, 
we briefly review relevant background from multilinear algebra. In 
\pref{sec:inv}, the invariance properties of a large class of objective functions, subsuming $\ploss$, are studied. In \pref{sec:sb}, the framework of symmetry breaking is introduced, and Puiseux series are constructed for various families of critical points using Newton polyhedra.

\section{Preliminaries}\label{sec:preliminaries}
A formal discussion of our results requires some familiarity with multilinear algebra. Parts of the exposition follow \cite{arjevani2021symmetry}.

A real tensor of order $n$ is an element of the tensor product of $n$ vector spaces $E_1 \otimes \ldots \otimes E_n$, where $E_i$ are real vector spaces. Upon choosing bases for each factor $E_i$, we may 
identify a tensor with a multi-dimensional 
array in $\RR^{d_1 \times \ldots \times d_n}$ with 
$d_i = \dim(E_i)$. We write 
$T_{i_1,\ldots,i_n}$ for the $(i_1,\ldots,i_n)$'th 
coordinate of $T$. Given vectors $\v_i \in \RR^{d_i}$, $i=1,\ldots,n$, 
we write $\v_1 \otimes \ldots \otimes \v_n$ for the 
outer product, namely the element 
in $\RR^{d_1} \otimes \ldots \otimes \RR^{d_n}$ such that  $(\v_1 \otimes 
\ldots \otimes \v_n)_{i_1,\ldots,i_n} = ({v_1})_ {i_1} \ldots 
(v_n)_{i_n}$, and $\v^{\otimes n} := \v \otimes \ldots \otimes 
\v$ for  $n$-times products. The Frobenius inner 
product of two tensors $T$ and $S$ of the 
same shape is defined by $\inner{ T, S}_F = \sum_{i_1,\ldots,i_n} 
T_{i_1\ldots i_n} S_{i_1 \ldots i_n}$. It is easy to 
verify that for the Frobenius inner product 
\begin{equation}\label{eq:frobenius_vector_product}
	\inner{ \v_1 \otimes
	\ldots \otimes \v_n, \w_1 \otimes \ldots \otimes 
	\w_n}_F = \prod_{i=1}^n \inner{\v_i, 
	\w_i}  
\end{equation}
for all $\v_i, \w_i \in \reals^{d_i}$.  In 
particular, $\inner{\v^{\otimes n}, \w^{\otimes 
n}}_F = \inner{\v,  \w}^n$.\\

A tensor $T \in (\RR^{d})^{\tens n}$ is \emph{symmetric} if it 
is invariant under permutation of indices, that is, if 
$T_{i_1,\ldots,i_n} = 
T_{i_{\sigma(1)},\ldots,i_{\sigma(n)}}$ for any
permutation $\sigma \in S_n$. The space of symmetric 
tensors of order $n$ on $\RR^d$, denoted by ${\rm \symp}^n(\RR^d)$, is 
$\binom{d+n-1}{n}$-dimensional and is isomorphic to the space of 
homogeneous polynomials of degree $n$ in $d$ variables. A natural 
isomorphism is given by the map 
\begin{align}\label{sp_id}
S\mapsto P_S,~~P_S(x) \defeq \inner{ S, 	x^ {\otimes r}}_F.
\end{align}

When $n=2$, this is the usual correspondence between 
symmetric matrices and quadratic forms (see \cite[Section 
3.1]{comon2008symmetric} for more details). The two norms considered in 
this 
work are induced from an inner product over ${\rm \symp}^n(\RR^d)$: the 
Frobenius inner 
product $\inner{\cdot,\cdot}_F$, and 
\begin{align}\label{data_ip}
	\inner{S, T}_{\cD} \defeq \mathbb{E}_{\x \sim \mathcal{D}} 
	[\inner{S, 
		\x^{\otimes 
			n}}_F \inner {T, 	\x^{\otimes n}}_F],
\end{align}	
where $\cD$ is a distribution on $\reals^d$ equivalent to Lebesgue 
measure, chosen so that the above expectation is finite for any $S,T\in 
\symp^n(\reals^d)$. 
\begin{fact}[\cite{arjevani2021symmetry}, Notation as above]\label{prop:data_ip} 
The bivariate function $\inner{\cdot, 
\cdot}_{\cD}:\symp^n(\reals^d) \times \symp^n(\reals^d)\to\reals$ defines an 
inner product on $\sympn$.
\end{fact}
Note that $\inner{\cdot,\cdot}_\cD$ is not positive definite over 
$\reals^{\tens d}$ and so the explicit restriction to $\sympn$ is 
necessary. Moreover, in terms of the coefficients of $P_S(x) = \sum 
a_{p_1,\ldots,p_d} x_1^{p_1}\ldots x_d^{p_d}$ and $P_T(x) = \sum 
b_{p_1,\ldots,p_d} x_1^{p_1}\ldots x_d^{p_d}$, the Frobenius inner product 
on $\symp^n(\RR^d)$ is given by
\begin{equation}\label{eq:frobenius_bombieri}
\inner{S, T}_F = \sum_{p_1+ \ldots + p_d = n} \binom{n}{p_1, \ldots,
p_d} a_{p_1,\ldots,p_d}	b_{p_1,\ldots,p_d},
\end{equation}
and the inner product $\inner{\cdot,\cdot}_\cD$ reads
\begin{equation}\label{eq:distributional_product}
	\inner{T, S}_{\mathcal D} = \sum_{\substack{p_1+ \ldots + p_d = n\\
			q_1+ \ldots + q_d = n}} \mathbb E_{\x \sim \mathcal{D}}[x_1^
	{p_1 + q_1}\ldots x_d^{p_d + q_d}] a_{p_1,\ldots,p_d} 
	b_{q_1,\ldots,q_d}.
\end{equation}
For the standard multivariate normal distribution $\cD= 
\cN(0,I_d)$, the above may be given explicitly,
\begin{align}
	&\inner{S, T}_{\cN} = \sum h_{\p+\q} a_
	{p_1,\ldots,p_d} b_{p_1,\ldots,p_d},
	\\
	&\mbox{ where } h_{\r} := 
	\begin{cases} 0 & \mbox{if } r_i \mbox{ is odd for some } i \in
		\{1,\ldots,d\},\\
		\prod_{i=1}^d {(r_i-1) !!} & \mbox{otherwise. }
	\end{cases}\nonumber
\end{align}
\begin{rem}
As noted in \cite{arjevani2021symmetry}, comparing 
\eqref{eq:frobenius_bombieri} and 
\eqref{eq:distributional_product}, it is seen that no data distribution induces 
the Frobenius inner product. 
Indeed, in~\eqref{eq:distributional_product}, whether a term 
$a_{p_1,\ldots,p_d}$ is multiplied by $b_{q_1,\ldots,q_d}$ depends only on 
${p_1 + q_1},\allowbreak\ldots,{p_d
	+ q_d}$, and so in particular if the coefficient of $a_ 
{p_1,\ldots,p_d}
b_{q_1,\ldots,q_d}$ is non-zero then, assuming, e.g., $p_1 > 0$, so is the 
coefficient corresponding to $a_{p_1-1,\ldots,p_d} b_{q_1+1,\ldots,q_d}$. 
In contrast, the coefficient of $a_ {p_1,\ldots,p_d} b_ 
{q_1,\ldots,q_d}$ in the Frobenius product is non-zero if and only if 
$(p_1,\ldots,p_d) = 
(q_1,\ldots,q_d)$. 
\end{rem}

A symmetric tensor $T \in \symp^n(\RR^d)$ is said to be of (real) rank 
one if $T = \alpha \w^{\otimes n}$ for some $\alpha \in \RR \setminus 
\{0\}$ and $\w \in \RR^d \setminus \{0\}$. A tensor is (real 
symmetric) rank $k$ if it can be written as a linear combination of $k$ 
rank-one tensors $T = \alpha_1 \w_1^{\otimes n} + \ldots + \alpha_k \w_k^{\otimes n}$, but not as a combination of $k-1$ rank-one tensors. For $n=2$, this definition agrees with the usual notion of rank of symmetric matrices. 

\section{Invariance Properties} \label{sec:inv}
Assume the target tensor $T$ defining (\ref{prob:opt}) is given by $T = \sum_{i=1}^{h}\beta_i \v_i^{\otimes n}$ for some $V\in M(h,d)$, $\bbeta\in\RR^h$ and $h\in\NN$. Making the dependence of $\ploss$ on $T$ explicit, it follows by \pref{prop:global} that the task of finding a rank $k$ decomposition may be equivalently stated as that of finding a global minimizer~for 
\begin{align} \label{prob:ext}
	\ploss(W,\balpha;V,\bbeta) = \nrm*{ \sum_{i=1}^{k}\alpha_i 
	\w_i^{\otimes n}- \sum_{i=1}^{h}\beta_i 
		\v_i^{\otimes n}}^2
\end{align}
under some tensor norm (in this work, the Frobenius and the cubic-Gaussian norm). Pairs of $(W,\balpha)$ with $\ploss(W,\balpha;V,\bbeta)=0$ therefore give different rank $k$ decompositions. The extended formulation of the loss function allows a natural description of invariance properties of 
$\ploss$ that follow by \emph{non-inherent }symmetries of $T$ (namely, ones occurring beyond those that are given by the mere definition of symmetric tensors).

\begin{rem}
If the rank of $T$ is exactly $k$, the above problem reduces to 
that of \textit{Waring decompositions} by  the identification 
given in (\ref{sp_id}) \cite{landsberg2011tensors}. 
\end{rem}

It will be convenient to adopt a general measure of similarity between vectors in $\reals^d$, here referred to as a \emph{kernel}. The use of kernels allows a succinct description of the assumptions used along the derivation of the invariance properties of $\ploss$, and yields general results, subsuming in particular settings considered in this work, with no additional cost. The choice to express invariance properties in terms of kernels may be motivated as follows. Observe that 
\begin{align*} 
&\ploss(W,\balpha;V,\bbeta) =
\inner{\sum_{i=1}^{d} 
		\alpha_i\w_i^{\otimes n} - \sum_{i=1}^{h} 
		\beta_i\v_i^{\otimes 
			n}, \sum_{i=1}^{d} 
		\alpha_i\w_i^{\otimes n} - \sum_{i=1}^{h} 
		\beta_i\v_i^{\otimes 
		n}}_F\nonumber	\\
%	&=  \sum_{i,j=1}^{k} \alpha_i\alpha_j\inner{\w_i^{\tens n}, 
%	\w_j^{\tens n}}_F
%	-2\sum_{i=1}^{k} \alpha_i\beta_j\sum_{j=1}^{h} 
%	\inner{\w_i^{\tens n}, \v_j^{\tens n}}_F
%	+ \sum_{i,j=1}^{h} \beta_i\beta_j\inner{\v_i^{\tens n}, 
%	\v_j^{\tens n}}_F\nonumber\\
	&=  \sum_{i,j=1}^{k} \alpha_i\alpha_j\inner{\w_i, \w_j}^n
	-2\sum_{i=1}^{k} \sum_{j=1}^{h} \alpha_i\beta_j\inner{\w_i, \v_j}^n
	+ \sum_{i,j=1}^{h} \beta_i\beta_j\inner{\v_i, \v_j}^n.\nonumber
\end{align*}
We may replace $\inner{\cdot, \cdot}^n$ by a 
general \emph{kernel} function $\ker:\reals^d\times\reals^d\to 
\reals$ satisfying $\ker(\w,\v) = \ker(\v,\w)$ and so define,
\begin{align} \label{prob:kernel}
\ploss_\ker&(W,\balpha;V,\bbeta) \defeq\\
&\sum_{i,j=1}^{k} \alpha_i\alpha_j\ker\prn{\w_i, \w_j}
-2\sum_{i=1}^{k} \sum_{j=1}^{h} \alpha_i\beta_j\ker\prn{\w_i, \v_j}
+ \sum_{i,j=1}^{h}\alpha_i\beta_j \ker\prn{\v_i, \v_j}.\nonumber
\end{align}
We note that in contrast to the conventional use of kernels in 
machine learning, in (\ref{prob:kernel}), $\w_i$ are \emph{optimization variables}. That is, both $\balpha$ and $W$ may be adjusted during the optimization process. Moreover, for the purpose of establishing the invariance properties of $\ploss$,  $\ker$ is not required to be positive definite. 

We shall be interested in a class of distribution-dependent kernels given by
\begin{align} \label{ker_cd}
\ker_\cD(\w,\v) = \bE_{\x \sim \cD} [\rho(\inner{\w,\x})\rho(\inner{\v,\x})],
\end{align}
with $\cD$ denoting a distribution over $\RR^d$ and $\rho:\reals\to\reals$ a measurable function. For example, the choice of $\rho(x)=\max\{x,0\}$, the ReLU activation function, and $\cD=\cN(0,I_d)$, is used in the study of two-layer ReLU networks 
\cite{ArjevaniField2019spurious}. 

\begin{rem} 
If the kernel function $\ker$ is positive definite, as ones considered in this work, $\ploss_\ker$ can be 
expressed in terms of inner products of the rows of $W$ and $V$ suitably mapped into the (unique) reproducing kernel Hilbert space corresponding to $\ker$ \cite{aronszajn1950theory}. 
\end{rem}

\begin{fact}[\cite{arjevani2021symmetry}, Notation as above]
	\label{prop:tensor_ip}
For the tensor norm (\ref{data_ip}), $\ploss \equiv \ploss_{\ker_\cD}$, the latter given by (\ref{prob:kernel}) with kernel function (\ref{ker_cd}) using $\rho(x)=x^n$.
\end{fact}
The data-dependent kernel corresponding to $\cD  = \cN(0,I_d)$ and $\rho(x)=x^3$ is referred to by the \emph{cubic-Gaussian 
kernel} \cite{arjevani2021symmetry}, and is given by
\begin{align}\label{kernel:cubic_gauss}
\ker_\cN(\w,\v) &= 
%\nu_3(\w,\v) &= 
	%	\bE_{\x}[\sigma(\inner{\w, \x})\sigma(\inner{\v,\x})] 
	%	&= \frac{9}{2}\prn*{ \frac{8}{6}\inner{\w,\v}^3 + 
		%	2\|\w\|^2\|\v\|^2\inner{\v,\w}}\\
	%	&=
	{6}\inner{\w,\v}_F^3 + 9\|\w\|_2^2\|\v\|_2^2\inner{\v,\w}_F.
\end{align}

\begin{definition}\label{def:perm_inv}
We say that a kernel $\ker$ is  \emph{permutation-invariant} if $\ker(\sigma 
\w,\sigma\v)\allowbreak=\ker(\w,\v)$ for all $\w,\v\in\reals^d,\sigma\in 
S_d$, with $S_d$ acting on $\RR^d$ by $(\sigma \w)_i = \w_{\sigma^{-1}(i)}$. 
\end{definition}
Since the Euclidean inner product $\inner{\cdot,\cdot}$ is $O(d)$-invariant and since the action of 
$S_d$ used in \pref{def:perm_inv} is orthogonal, $\ker(\w,\v) = 
\inner{\w, 
\v}^n$ is permutation-invariant for any $n\in\bN$. In a similar vein, any 
kernel expressed in terms of $\inner{\w,\v}$, such as  $\ker_{\cN}$, is 
permutation-invariant.

By a direct extension of the action of $S_d$, we let a pair of permutations 
$\sigma = (\sigma_1,\sigma_2)\in S_{k} \times 
S_{d}$ act on a matrix $A = \left[A_{ij}\right]
\in M(k,d)$ by
\begin{equation} \label{sksd_action}
	\left(\sigma A\right)_{ij} \defeq
	A_{\sigma_1^{-1}(i),\sigma_2^{-1}(j)}.
\end{equation}
Thus, $\sigma_1$ acts by permuting the rows of $A$, and $\sigma_2$ by 
permuting the columns of $A$. The largest subgroup of $S_{k} \times 
S_{d}$ fixing $A\in M(k,d)$, namely, the isotropy subgroup of $A$, is denoted by $(S_k\times S_d)_A$. Likewise, $\prn{S_d}_{\x}$ is 
the largest subgroup of $S_d$ fixing $\x\in\reals^d$. For example, the isotropy 
subgroup of $I_d\in M(d,d)$ is $\Delta S_{d}$, $\Delta$ mapping
a subgroup $G\subseteq S_d$ to the diagonal subgroup $\Delta G 
\defeq \{(g,g) : g \in G\}$. The isotropy group of $\cI_{d,1}$, 
$\cI_{d_1,d_2}$ generally denoting the $d_1\times d_2$-matrix with all entries equal to 1, is $S_d$. We denote the projection of a subgroup $H\subseteq S_{d_1}  \times S_{d_2}$ onto its $i$'th component by $\pi_i(H)$. In 
particular, $\pi_1(\Delta G) = \pi_2(\Delta G) = G$ for any subgroup~$G\subseteq S_d$.

\begin{prop}[Notation and assumption as above]\label{prop:inv}
Fix $V\in M(h,d)$, $\bbeta\in\reals^h$, $\balpha\in\reals^k$ and define 
$\isov_{\bbeta,V} \defeq \Pi_2((S_h\times S_d)_{V}\cap ((S_h)_{\bbeta}\times 
S_d) 
)$. If a kernel $\ker$ is permutation-invariant then $\ploss_\ker$ is $(S_k)_{\balpha}\times 
\isov_{\bbeta,V}$-invariant with respect to $W$.
That is, $\ploss_\ker(\sigma W, \balpha; V,\bbeta) = \ploss_\ker(W, \balpha; 
V,\bbeta)$ for all $\sigma\in (S_k)_{\balpha}\times \isov_{\bbeta,V}$.
\end{prop}
\pref{prop:inv} generalizes the invariance results given in 
\cite{arjevani2021symmetry,ArjevaniField2019spurious}, giving  an explicit description of the invariance properties of permutation-invariant kernels in terms of the symmetry of $V$, $\balpha$ and $\bbeta$. It is seen that the larger the isotropy of the target matrix $V$, the richer the invariance properties of $\ploss_\ker$. 
For example, if $\balpha = \bbeta = \cI_{d,1}$  and $V \in M(d,d)$ is a 
circulant matrix then $(S_d\times S_d)_{V}\supseteq \Delta \ZZ_d$, with $\ZZ_d 
\defeq  \ZZ/d\ZZ$ denoting the group of all cyclic permutations, rendering
$\ploss_\ker$ $S_k\times \ZZ_d$-invariant. A similar argument gives the 
following.
\begin{cor} [Notation and assumptions as above]\label{corr:ploss_inv}
If $V=I_d\in M(d,d)$, $\bbeta = \cI_{d,1}$ (hence defining $T_\e$) and $\balpha = \cI_{k,1}$, then optimization problem (\ref{prob:opt})  is $S_k\times S_d$-invariant under the Frobenius as well as the cubic-Gaussian tensor norms.
\end{cor}
The corollary follows by a direct application of \pref{prop:inv}, 
noting  that 
$\isov_{\bbeta,V} = \pi_2((S_d\times S_d)_{V}\cap ((S_d)_{\bbeta}\times S_d)) = 
\pi_2(\Delta 
S_d \cap (S_d \times S_d)) = \pi_2(\Delta S_d)=S_d$.

Given $G\subseteq S_k\times S_d$, 
the gradient vector field of a $G$-invariant differentiable function 
is $G$-equivariant (an inner product on $M(k,d)$ is tacitly 
understood) and so maps $M(k, d)^G$ to itself. By \pref{prop:inv}, 
$\ploss_\ker$ is $(S_k)_{\balpha}\times \isov_{\bbeta,V}$-invariant with 
respect to~$W$. In particular, $\nabla_W \ploss_\ker$ restricted to 
$M \defeq M(k,d)^{(S_k)_{\balpha}\times \isov_{\bbeta,V}}$  satisfies
\begin{align}\label{grad_equiv}
(\nabla_W \ploss_\ker)|_{M}(\sigma W, \balpha; V,\bbeta) =
\sigma(\nabla_W \ploss_\ker)|_{M}(W, \balpha; V,\bbeta),
\end{align}
for all $\sigma\in (S_k)_{\balpha}\times \isov_{\bbeta,V}$.

\section{Families of symmetry breaking critical points}\label{sec:sb}
By \pref{corr:ploss_inv}, if $k=d $ and if the target tensor used is $T_\e$ then $\ploss$ is $S_d\times S_d$-invariant. As indicated earlier, the isotropy groups of critical points of $\ploss$ found by numerical methods are typically proper subgroups of $S_d\times S_d$ and so are symmetry breaking \cite{arjevani2021symmetry}. Below, we outline a framework developed in \cite{arjevani2023hidden,arjevanifield2020hessian,ArjevaniField2021,ArjevaniField2019spurious,ArjevaniField2020} for studying phenomena of symmetry breaking of this  nature. Our presentation makes extensive use of a Newton polyhedron argument which is new to the framework and allows a systematic investigation of all critical points occurring in a given fixed points space. The methods are first illustrated with reference to the $S_d\times S_d$-isotropy case under the Frobenius norm, and are later described for general isotropy groups and tensor norms. Throughout, we assume that the target tensor is $T_\e$ (but see \pref{sec:generalities}).

Fix $d\in\NN$ and let $\fC(d)\in 
M(d,d)^{S_d\times S_d}$ be a critical point of $\ploss_F$. 
Since $M(d,d)^{S_d\times S_d} = \{\xi \cI_{d,d}~|~\xi \in \reals\}$,  $\fC(d) = \xi(d) 
\cI_{d,d}$ for some real scalar $\xi(d)\in\RR$. The condition $	
\nabla \ploss_F (\fC(d)) = 0$, omitting the dependence of $\ploss_F$ on $V, \balpha$ and $\bbeta$ for convenience, reads 
\begin{align*}
	\nabla \ploss_F (\xi(d) \cI_{d,d}) = 0,
\end{align*}
giving a single polynomial equation in $\xi$ and $d$,
\begin{align}\label{eqn:single}
\xi^5d^3-\xi^2= 0.
\end{align}
Thus, for any fixed $d$, there are two (real) critical points with isotropy 
$S_d\times S_d$: $\fC_0(d) = 0_{d,d}$ and $\fC_1(d) = d^{-1}\cI_{d,d}$. 
Moreover,
\begin{align}
\ploss_F(\xi \cI_{d,d}) &= \xi^6d^5 - 2\xi^3d^2+d,
\end{align}
and so the respective loss is 
$\ploss_F(\fC_0(d)) = d$ and $\ploss_F(\fC_1(d)) = d - \frac{1}{ d}$. The loss, in both cases, is strictly positive for $d\ge2$ and so in particular $\fC_1$ cannot be a global minimum. To determine whether $\fC_1$ defines a local minimum, a local maximum or a saddle, we compute the Hessian spectrum. The symmetry of $\fC_1$ allows a considerable simplification in the computation of the Hessian spectrum using results from the representation theory of the symmetric group \cite{fulton1991representation}. The analysis proceeds by considering four irreducible $S_d$-representations (isomorphism classes of $S_d$-invariant linear spaces containing no proper $S_d$-invariant subspace but the trivial):

\begin{itemize}[leftmargin=*]
	\item The trivial representation $\mathfrak{t}_d$ of degree 	1.
	\item The standard representation $\mathfrak{s}_d$ of $S_d$ of degree $d-1$.
	\item The exterior square representation $\mathfrak{x}_d=\wedge^2 
	\mathfrak{s}_d$ of degree $\frac{(d-1)(d-2)}{2}$.
	\item The representation $\mathfrak{y}_d$  associated to the 
	partition $(d-2,2)$ of degree~$\frac{d(d-3)}{2}$.
\end{itemize}

The Hessian spectrum is computed and given separately for each irreducible 
representation or tensor products thereof (see \pref{table:sdTsdevs} and 
\pref{table:normal_sdTsdevs}). We refer to \cite[Section XII]{golubitsky2012singularities} for more details on the approach.

In the general case, one considers sequences of subgroups 	
$(G_d)_{d\in\NN}$ of $S_k\times S_d$  for which  $\dim(M(k, d)^{G_d})$ stabilizes for sufficiently large $d$. To simplify exposition, assume $\dim(M(k, d))$ is fixed and equal to $N\in\NN$ for all $d\in\NN$. The fixed point space bears a natural parameterization given by a linear isomorphism 
$\Xi \defeq \Xi^{G_d}:\RR^N\to M(k, d)^{G_d}$. Families of critical points in $M(k,d)^{G_d}$ may therefore be expressed by $\fC(d) = \Xi(\fc(d))$, with $\fc(d)\in\RR^N$ satisfying a system of $N$ polynomial equations, 
\begin{align}\label{grad_eqs}
\Xi^{-1}\circ\nabla \ploss|_{M(k, d)^{G_d}} \circ 
\Xi(\fc(d)) = 0.
\end{align}
Note that pulling $\nabla \ploss|_{M(k, d)^{G_d}}$ back to $\RR^N$ via $\Xi^{-1}$ is possible by the $G_d$-equivariance of the gradient vector field $\nabla \ploss$, see (\ref{grad_equiv}). The single polynomial equation given in the $S_d\times S_d$-case, namely \pref{eqn:single}, now becomes a {system }of polynomial equations, and so expressing $\fC(d)$ in terms of simple formulae in $d$ is generally hard. However, for the purpose of 
characterizing the stability of families of critical points it suffices, in many cases, to estimate the Hessian spectrum to a certain order. An exact computation of $\fC(d)$ is not required---estimates to a limited order suffice. Here estimates are obtained by considering partial sums of Puiseux series constructed by repeating a Newton polytope argument (details appear in \pref{sec:rank_for}). In \pref{sec:puiseux}, we provide a detailed construction of a family of (complex) critical points in Puiseux series in $d^{\frac{-1}{4}}$ based on a partial (fractional) sum so obtained. General results on the existence of Puiseux series solutions are given in \cite{mcdonald2002fractional}. To our applications, the CSL, applied to the semialgebraic solution set (\ref{grad_eqs}), suffices \cite{arjevani2024symmetry}. A more direct approach is given by applying the CSL to the solution set transformed by (e.g., \cite{chen2014milnor})
\begin{align}
\bxi \defeq
(\xi_1,\dots,\xi_N) &\mapsto 
\prn*{
\frac{\xi_1}{\sqrt{1+\|\bxi\|^2}}
,\dots,
\frac{\xi_N}{\sqrt{1+\|\bxi\|^2}},
\frac{1}{\sqrt{1+\|\bxi\|^2}}},
\end{align}
see \cite{arjevani2024symmetry} (also addressing the dimensionality of strata of the solution~set).

\subsection{Critical points under the Frobenius norm} \label{sec:rank_for}
For the target tensor $T_\e$, the loss function 
$\ploss_F$ takes the form
\begin{align}
\varphi(W) \defeq \nrm*{\sum_{i=1}\w_i^{\tens 3} - \sum_{i=1}\e_i^{\tens 
3}}_F^2.
\end{align}

\subsubsection{Isotropy $\Delta S_d$.} Critical points of $S_d\times 
S_d$-isotropy are given in the previous section and so we begin by 
considering the $\Delta S_d$-case.
The fixed point space corresponding to $\Delta S_d$ is 
\begin{align*}
M(d,d)^{\Delta S_d}\defeq \{\xi_1 I_d + \xi_2 (\cI_{d,d}-I_d)~|~\xi_i\in\RR\}.	
\end{align*}
By the preceding section, we seek a solution of the form
\begin{align}\label{tau_lot}
\xi_i(d) = A_{i} d^{\tau_i}+ \text{terms with lower exponents}, i=1,2,
\end{align}
with $A_{i}\neq 0$ (vanishing $A_i$ will be discussed later) and 
$\tau_i\in\QQ$ (note that exponents are assumed decreasing) for a given family of critical points. The critical point equations corresponding to the restriction of $\grad \varphi$ to $M(d,d)^{\Delta S_d}$ (divided by 2) are
\begin{align}\label{deltasd_cubic}
0 &= 3 
\xi_{2}^{5} d^{3}  + d^{2} \left(15 \xi_{1} 
\xi_{2}^{4} - 15 \xi_{2}^{5}\right) + d \left(6 \xi_{1}^{3} 
\xi_{2}^{2} + 12 \xi_{1}^{2} \xi_{2}^{3} - 42 \xi_{1} 
\xi_{2}^{4} + 24 \xi_{2}^{5}\right)\nonumber\\
&+3 \xi_{1}^{5} - 6 \xi_{1}^{3} \xi_{2}^{2} - 12 \xi_{1}^{2} 
\xi_{2}^{3} - 3 \xi_{1}^{2} + 27 \xi_{1} \xi_{2}^{4} - 12 
\xi_{2}^{5},\nonumber \\
0&= 3 \xi_{2}^{5} d^{3}  + d^{2}
\left(15 \xi_{1} 
\xi_{2}^{4} - 15 \xi_{2}^{5}\right) + d \left(30 \xi_{1}^{2} 
\xi_{2}^{3} - 60 \xi_{1} \xi_{2}^{4} + 30 \xi_{2}^{5}\right)\nonumber\\
&+ 3 \xi_{1}^{4} \xi_{2} + 12 \xi_{1}^{3} \xi_{2}^{2} - 54 
\xi_{1}^{2} \xi_{2}^{3} + 60 \xi_{1} \xi_{2}^{4} - 21 
\xi_{2}^{5} - 3 \xi_{2}^{2}.
\end{align}
When $\xi_i$ in the above are substituted by (\ref{tau_lot}), every term of maximal total degree must appear at least twice so that a cancellation by a different term of the same (total) degree is possible by a suitable choice of $A_i$, see e.g., \cite{sturmfels2002solving}. Equivalently, 
arguments 
maximizing 
\begin{align}\label{pl_cond}
\max\left(2 \tau_{1}, 5 \tau_{1}, 5 \tau_{2} + 3, \tau_{1} 
+ 4 \tau_{2} + 2, 2 \tau_{1} + 3 \tau_{2} + 1, 3 \tau_{1} + 2 \tau_{2} + 
1\right)
\end{align}
must be attained at least twice, as well as 
arguments maximizing  
\begin{align*}
&\max\left(2 \tau_{2}, 3 \tau_{1} + 2 \tau_{2}, 4 \tau_{1} 
+ \tau_{2},5 \tau_{2} + 3, \tau_{1} + 4 \tau_{2} + 2, 2 \tau_{1} + 3 \tau_{2} 
+ 1\right).
\end{align*} 
The only possibilities are $(\tau_1,\tau_2) = 
(0,-3/4)$ and $(\tau_1,\tau_2) = (-1,-1)$. Since $A_i$ are chosen so that terms of maximal total degree cancel out,  for $(\tau_1,\tau_2) = (-1,-1)$, $(A_1, A_2) = (\pm1, 1)$, 
defining two families of minima: $\fC_1$ of isotropy  $S_d\times S_d$
encountered earlier and $\fC_2(d) = (1/d)(\cI_{d,d} - 2I_d)+o(d)$ of 
isotropy $\Delta S_d$. For $(\tau_1,\tau_2) = (0, -3/4)$, $(A_1, 
A_2) = (1, \sqrt[4]{-1})$, where $\sqrt[4]{-1}$ denotes any fourth root of $-1$, 
none of which is real. The families of critical 
points $\fC_0(d) \defeq 0_{d,d}$
and $\fC_I(d) \defeq I_{d,d}$ are obtained by letting either $A_i$ or 
both vanish. Additional terms of the Puiseux series are computed by 
repeating the process.

\begin{rem}\label{rem:complex}
Depending on the application, tensor decomposition problems may be studied over 
the complex numbers, e.g., 
\cite{comon2002,cartwright2013number,chiantini2017generic,brachat2010symmetric}.
 The symmetry breaking framework carries over to this setting without change. 
For example, a family of complex $\Delta S_d$-critical points corresponding to 
the pair of exponents $(\tau_1,\tau_2) = (0,-3/4)$ is given by 
\begin{align*}
(\fc_{\bC})_{1}(d) &= 1 - \frac{2 i}{3 d^{\frac{ 1}{2}}} + \frac{\sqrt{2} \left(1 + i\right)}{d^{\frac{3}{4}}} - \frac{77}{36 d} + o(d),\\
(\fc_{\bC})_{2}(d) &=  \frac{\sqrt{2} \left(1 + i\right)}{2 d^{\frac{3}{4}}}- \frac{5}{4 d} + o(d).
\end{align*}
Note that here Puiseux series are given in terms of $d^{-\frac{1}{4}}$.
\end{rem}

\begin{rem}
By B\'{e}zout's theorem the number of solutions of (\ref{deltasd_cubic}) is generically 25 (the product of the total degree of the equations), exceeding the actual number of solutions. The more refined bound given by Bernstein's theorem \cite{bernshtein1975number} is still lose. In both cases, this results from additional structure on account of symmetry, see \cite{arjevani2024symmetry}.
\end{rem}

\subsubsection{Isotropy groups with two or more factors.} 
For critical points of isotropy $\Delta (S_{d-1}\times S_1)$, the 
analog for (\ref{pl_cond}) comprises five piecewise-linear conditions (with seventy two terms). %The system, as well as its thirteen 
%solutions, is given in \pref{sec:solution_dsd}. 
A $\Delta (S_{d-1}\times S_1)$-family of critical points, denoted $\fC_3 = \Xi^{\Delta 
	(S_{d-1}\times S_1)}(\fc_3)$, is given by 
\begin{align*}
(\fc_3)_{1}(d) &= - \frac{1}{d} - \frac{13}{3 d^{2}} - \frac{77}{9 d^{3}} + 
\frac{6421}{81 d^{4}} + o(d^{-4}),&(\fc_3)_4(d) &= o(d^{-4}),\\
(\fc_3)_2(d) &=\frac{1}{d} + \frac{13}{3 d^{2}} + \frac{149}{9 d^{3}} + 
\frac{2867}{81 d^{4}} + o(d^{-4}),&(\fc_3)_5(d) &= 1 + o(d^{-4}),\\
(\fc_3)_3(d) &= o(d^{-4}).
\end{align*}
In some cases, the entries of a family of critical points do not depend 
on $d$. One such example is $\fC_4 = \Xi^{\Delta (S_{d-2}\times 
S_2)}(\fc_4)$  
where
\begin{align}
(\fc_4)_1(d) =1, ~(\fc_4)_5(d) = 1/2,~(\fc_4)_6(d) = 1/2
\end{align}
(only non-zero entries are described). Additional examples are given by $\fC_0$,  $\fC_I$ and $\fC_{5,t}$, defined in \pref{sec:global}.

\subsection{Critical points for the Cubic-Gaussian norm} 
\label{sec:cubic_gauss}
As earlier, we assume that the target tensor is $T_\e$. The loss function $\ploss_\cN$  takes the form
\begin{align}
	\psi(W) \defeq \nrm*{\sum_{i=1}\w_i^{\tens 3} - \sum_{i=1}\e_i^{\tens 
			3}}_\cN^2.
\end{align}
\subsubsection{Isotropy $S_d\times S_d$.}
We begin by considering critical points in $M(d,d)^{S_d\times S_d} = \{\xi_1 \cI_{d,d}~|~\xi \in\reals\}$. We have, 
\begin{align*}
\psi(\xi I) &= 15 \xi^{6} d^{5} - 18 \xi^{3} \left(d^{3} 
+ \frac{2 d^{2}}{3}\right) + 15 d,\\
(\nabla \psi)|_{M(d,d)^{S_d\times S_d}} &= 90 \xi^{5} d^{3} - 54 
\xi^{2} \left(d + \frac{2}{3}\right).
\end{align*}
There are therefore two families of critical points corresponding to isotropy $S_d\times S_d$, $\fD_0(d) = 0_{d,d}$ and 
\begin{align*}
\fD_1(d) = \sqrt[3]{ \frac{3}{5d^3}  \left(d + \frac{2}{3}\right)	} \cdot \cI_{d,d}.
\end{align*}
The loss at $\fD_0$ is $\varphi(\fD_0(d)) = 15d$. The Hessian at $\fD_0$ vanishes (but similarly to $\fC_0$, the tensor of third-order partial derivatives does not vanish). The loss at $\cD_1$ is $\psi(\cD_1) = \frac{48 d}{5} - \frac{36}{5} - \frac{12}{5 d}$. The Hessian spectrum is given in \pref{table:normal_sdTsdevs}.

\subsubsection{Diagonal Isotropy groups.}
For isotropy $\Delta S_d$, the system of piecewise-linear conditions analogous 
to (\ref{pl_cond}) implies that a critical point given by 
(\ref{tau_lot}) must have $(\tau_1, \tau_2) = (0,-1)$, $(\tau_1, \tau_2) = 
(-1/6,-2/3)$ or $(\tau_1, \tau_2) = (-2/3,-2/3)$.
The first and the second case imply  $A_1A_2=0$, and are therefore 
disregarded (recall that $A_i$ are assumed to be non-zero). The third case gives $\fD_1$.  

The family of critical points for which $A_2=0$ is $\fD_{I}(d)  = I_d$. Clearly, $\psi(\fD_{I}) = 0$. The Hessian spectrum is given in \pref{table:normal_sdTsdevs}. A family of $\Delta (S_{d-1}\times S_1)$-critical points $\fD_{2} = \Xi^{\Delta (S_{d-1}\times S_1)}(\fd_2)$ giving spurious minima is defined, modulo 
$o(d^{-4})$-terms, by 
\begin{align*}
(\fd_{2})_1(d) &=1 - \frac{5}{3 d^{3}} + \frac{9}{d^{4}}, &(\fd_{2})_4(d) = 0,\\ 
(\fd_{2})_2(d) &=- \frac{1}{d^{3}} + \frac{7}{d^{4}}, &(\fd_{2})_5(d) = 0,\\
(\fd_{2})_3(d) &= 	\frac{1}{d} - \frac{1}{d^{2}} + \frac{2}{d^{3}} + \frac{10}{3 d^{4}}.
\end{align*}

\section{Concluding remarks}\label{sec:concluding}
Standard approaches \cite{carroll1970analysis,harshman70,paatero1997least,tomasiuse,acar2011scalable,kolda2015numerical} for decomposing real symmetric tensors often fail to return a valid decomposition \cite{kolda2009tensor,acar2011scalable,kolda2015numerical,tomasi2006comparison} of nonconvexity of the associated optimization landscape. Our focus in this article has been on studying in some depth obstructions to local optimization methods using techniques pertaining to symmetry breaking~\cite{arjevani2024symmetry}. The complex array of saddles and minima so revealed exhibits various phenomena that influence the efficacy of local optimization methods. Some phenomena are desirable; for instance, saddles with higher objective values tending to exhibit a greater number of descent directions, as well as the emergence of such directions under over-parametrization (see \emph{op. cit.} and \cite{ArjevaniField2022annihilation}). Others are awkward,  as third-order saddles existing in the optimization landscape, which can unduly slow down the optimization process or even cause an optimization method to erroneously declare convergence to a local minimum. Certain behaviors are common to different tensor norms, such as highly skewed Hessian spectra, while others, like the growth rates of Hessian eigenvalues, differ. These phenomena, along with others discussed in detail in \pref{sec:main_results}, help clarify key questions foundational to theories that aim to characterize the decomposability of symmetric tensors using local optimization methods. We conclude by noting that, although our emphasis has been on orthonormally decomposable tensors, this assumption is not strictly necessary. Variations, beyond those discussed in \pref{sec:main_results}, include, for instance, further emphasis on collinearity factors \cite{li1997numerical}, over-factoring \cite{tomasi2006comparison}, and tensors with generic rank \cite{hirschowitz1995polynomial}.

In a broader context, the applicability of the symmetry-breaking framework stems from the presence of critical points whose isotropy groups are large. A formal mathematical mechanism accounting for the emergence of such symmetry-breaking points is developed in~\cite{arjevani2024symmetry}.

\section*{Acknowledgements}
We thank Michael Field, Yuval Peled and Avi Wigderson for several helpful and insightful discussions. The research was supported by the Israel Science Foundation (grant No. 724/22).

\bibliographystyle{ieeetr}
\bibliography{bib}

\appendix

% !TEX root = tensor_sym_crit_r.tex

\newpage

\section{Generalities} \label{sec:generalities}
Following \pref{sec:inv}, we consider a formulation of the objective function 
$\ploss_\ker:M(k,d)\to\RR$ (having dependence on other problem parameters implicit) given by a kernel $\ker:\RR^d\times \RR^d \to \RR$ satisfying $\ker(\w,\v) = 
\ker(\v,\w)$,
\begin{align}\label{eq:loss}
\ploss&(W) =
\sum_{i,j=1}^{k} \ker(\w_i, \w_j)
-2\sum_{i=1}^{k}  \sum_{j=1}^{h} \ker(\w_i, \v_j)
+ \sum_{i,j=1}^{h} \ker(\v_i, \v_j),
\end{align}
where $\w_i$ (resp. $\v_i$) is a column vector denoting the $i$'th row of $W\in M(k,d)$ (resp. 
$V\in M(h,d)$). The relevant literature considers different conventions 
for arranging first-, second- and higher-order derivatives of a real-valued 
function defined over matrices, see 
e.g.,  \cite{magnus2019matrix}. Here we find it convenient to identify 
$M(k,d)$ with $\RR^{kd}$ and regard $W\in M(k,d)$ as a column vector having 
 $kd$ entries given by the rows of $W$ transposed and stacked on top 
of one another. Thus, the differential $D\ploss_\ker$ maps $\RR^{kd}$  to 
$L(\RR^{kd},\reals)$, 
where $L(E,F)$ generally denotes the space of all bounded linear operators from 
a  vector space $E$ to a vector space $F$. The gradient $\nabla 
\ploss_\ker:\RR^{kd}\to \RR^{kd}$ (defined by the usual identification) is given~by
\newcommand{\vh}{{h}}
\begin{align} \label{eq:grad}
	\nabla \ploss(W) &= 2\sum_{i=1}^{k} \e_i\otimes \prn*{\sum_{j=1}^k  
		\kappa_\w
		(\w_i,\w_j)	- \sum_{j=1}^\vh \kappa_\w (\w_i,\v_j)} ,
\end{align}
where $\e_i$ denotes the $i$-th unit vector and $\kappa_\w:\RR^d\times \RR^d\to 
\RR^d$ a column vector containing the partial derivatives of $\kappa$ with 
respect 
to the entries of the first 
argument. The Hessian  $\nabla \ploss^2_\ker:\RR^{kd}\to 
L(\RR^{kd}, L(\RR^{kd},\RR)) \cong M(\RR^{kd},\RR^{kd})$ is given by 
\begin{align}\label{eq:hess}
\nabla^2 \ploss(W) 
&= 2\sum_{i,j=1}^{k} \e_i \e_j^\top\otimes \kappa_{\w,\v} (\w_i,\w_j)	
\\
&+ 2\sum_{i=1}^{k} \e_i \e_i^\top\otimes 
\prn*{\sum_{j=1}^k \kappa_{\w,\w} (\w_i,\w_j) - \sum_{j=1}^\vh 
	\kappa_{\w,\w} 
	(\w_i,\v_j)} ,\nonumber
\end{align}
where $\kappa_{\w, \w}$ (resp. $\kappa_{\w, \v}$) is the  $d\times d$ Jacobian 
matrix of 
$\kappa_\w$ with respect to $\w$ (resp. $\v$).

If $\ker(\w,\v)=\inner{\w,\v}^n$ and $h=d$, Equation \pref{eq:loss} reads
\begin{align}\label{eqn:loss}
\ploss&(W) =
\sum_{i,j=1}^{k} \inner{\w_i, \w_j}^n
-2\sum_{i=1}^{k}  \sum_{j=1}^{d} \inner{\w_i, \v_j}^n
+ \sum_{i,j=1}^{d} \inner{\v_i, \v_j}^n.
\end{align}
If $\v_1,\dots,\v_d$ are orthonormal, $V$ denotes the matrix whose rows 
are $\v_i$ and $\hat{\w}_i \defeq ({\w}_i V)^\top$, then
\begin{align} \label{eqn:obj}
\hat{\ploss}&(W) \defeq \ploss(WV) 
= \sum_{i,j=1}^{k} \inner{\hat{\w}_i, \hat{\w}_j}^n
-2\sum_{i=1}^{k}  \sum_{j=1}^{d} \inner{\hat{\w}_i, \e_j}^n
+ \sum_{i,j=1}^{d} \inner{\e_i, \e_j}^n.
\end{align}
Thus, the optimization landscape of ${\ploss}$ and $\hat{\ploss}$, the latter corresponding to the target tensor $T_\e$,  are identical modulo orthonormal transformation of the domain (a similar argument shows that this applies to Gauss-cubic norm as well). It is therefore no loss of generality to focus on the optimization landscape of $\hat{\ploss}$. Equation \pref{eq:grad} now reads
\begin{align} \label{eq:grad_red}
\nabla \ploss(W) &= 2n\sum_{i=1}^{k} \e_i\otimes \prn*{\sum_{j=1}^k  
\inner{\w_i, \w_j}^{n-1}\w_j - \sum_{j=1}^d \inner{\w_i, 
\e_j}^{n-1}\e_j}.
\end{align}

It follows by the above shows that eigenvectors of the target tensor $T$ 
\cite{lim2005singular,qi2005eigenvalues,robeva2016orthogonal} give
critical points. Indeed, a vector $\w\in\RR^d$ is an eigenvector of $T$ if $T 
\w^{n-1}=\lambda \w$ for some $\lambda\in\RR$, where $T 
\w^{n-1}$ generally denotes a  $d$-dimensional vector given by $(T 
\w^{n-1})_i \defeq 
\sum^d_{i_2,\dots,i_n = 1} T_{i, i_2,\dots,i_n} w_{i_2}\cdots 
w_{i_n}$. 
Given an eigenpair $(\w,\lambda)$, define a matrix $W_t$ whose first row is 
$t\w, 
t\in\RR$, and is otherwise zero. We have,
\begin{align}\label{eq:tscalar}
\nabla \ploss(W_t) 
%&= 2n \e_1\otimes \prn*{
%t\inner{t\w, t\w}^{n-1}\w	- tT\w^{n-1}}\\
&= 2n \e_1\otimes \prn*{
t^{2n-1}\inner{\w, \w}^{n-1}\w	- t^{n-1}\lambda \w}.
\end{align}
Thus, $W_t$ defines a critical point for 
$t=(\lambda/\inner{\w,\w}^{n-1})^{1/n}$. 
Generalizing 
the latter argument, if $r$ orthogonal eigenvectors $\w_1,\dots,\w_r$ are 
given, $t_i$ may 
be set individually as in \pref{eq:tscalar}. 
Forming a weight matrix $W_{(r)}$ whose first $r$ rows are 
$t_i\w_1,\dots,t_r\w_r$  
and is otherwise zero, we have,
\begin{align} \label{eqn:grad}
\nabla \ploss(W_{(r)}) &= 2n\sum_{i=1}^{r} \e_i\otimes \prn*{
	t_i^{2n-1}\inner{\w_i, \w_i}^{n-1}\w_i	- t_i^{n-1}\lambda \w_i} = 0.
\end{align}
Of course, families of critical points whose weight matrices are not 
orthonormal, such as 
$\fC_2$, do not arise by this mechanism. 

\begin{rem}
An equivalent way of defining an eigenvector follows using
$f_T:\RR^d\to\RR:w \mapsto \inner{T,w^d}_F$. A simple argument shows that

$(\w,\lambda)$ is an 
eigenpair iff it satisfies
\begin{align*}
\nabla f_T(w) = d\lambda w.
\end{align*}
It is possible to replace the Frobenius inner product in definition of $f_T$ with a different inner product and obtain similar results.
\end{rem}

The expressions of the Hessian \pref{eq:hess} are given by,  
\begin{align}
\nabla^2& \ploss(W) 
= 2\sum_{i,j=1}^{k} \e_i \e_j^\top\otimes 
\prn{n\inner{\w_i,\w_j}^{n-1}I_d + 
n(n-1)\inner{\w_i,\w_j}^{n-2}\w_j\w_i^\top}
\\
&+ 2\sum_{i=1}^{k} \e_i \e_i^\top\otimes 
\prn*{\sum_{j=1}^k n(n-1)\inner{\w_i,\w_j}^{n-2}\w_j\w_j^\top - \sum_{j=1}^d 
	n(n-1)\inner{\w_i,\v_j}^{n-2}\v_j\v_j^\top}. \nonumber
\end{align}
Analytic expressions for the Gauss-cubic norm are derived 
similarly.

\section{Proofs}
\subsection{Proofs for \pref{sec:global}}
\subsubsection{\pfref{prop:global}}
For any mode of the target tensor $T_\e = \sum_{i=1}^d \e_i^{\otimes 3}$, 
the matrix containing the respective columns is simply the identity 
matrix $I_d$. The 
rank 
of the identity matrix $I_d$ is $d$ and so trivially every $d$ columns are 
independent. 
Thus, by Kruskal's identifiability result \cite{kruskal1977three}, 
$T_\e$ admits a unique decomposition modulo permutation and scaling of the 
columns if $3d\ge 2d + 2$, i.e., if $d \ge 2$.

We also provide the details for $T = \sum_{i=1}^d \e_i^{\otimes n},n\ge 4$. 
For even $n$, we assume that  $\balpha$ is adjusted during the optimization process. If $(W, \balpha)$ is a global minimizer, then $\ploss(W) = 0$ and so independently of the underlying tensor norm, $T = \sum_{i=1}^d \alpha_i 
\w_i^{\otimes n}$. Rearranging indices if needed, we may assume 
$\w_i=c_i\e_i,~c_i\in\RR$ for $i\in[d]$, giving
\begin{align*}
\alpha_j c_j^n = \prn*{\sum_{i=1}^d \alpha_i c_i^n\e_i^{\otimes 
n}}_{j,\dots,j} = 
\prn*{\sum_{i=1}^d \alpha_i \w_i^{\otimes n}}_{j,\dots,j}
= T_{j,\dots,j} = 1.
\end{align*}
Therefore, if $n$ is odd, $c_j = \alpha_j^{-1/n}$. If $n$ is even, necessarily 
$\alpha_j \ge 0 $ and $c_j = \pm \alpha_j^{-1/n}$. (For $n=3$, we 
assumed $\alpha_j=1$ and so $c_j=1$.)

\subsection{\pfref{prop:proper}}
The proposition follows by a direct computation.
First, we show that $\fC_{5,t}(d) \defeq I_{d-2} \oplus \begin{bsmallmatrix} 
(1-t^3)^{1/3}&	0\\ t & 0\end{bsmallmatrix}$ is a critical point for 
every $t\in\RR$ and $d\in \NN$. As usual, we let $\w_i$ denote the rows of 
$\fC_{5,t}$. Referring to \pref{eqn:grad}, since the first $d-2$ rows are 
orthogonal to other rows, we have
\begin{align*}
\nabla \ploss(\fC_{5,t}(d)) &= 
2n\sum_{i=d-1}^{d} \e_i\otimes 
\prn*{\sum_{i=d-1}^d
	\inner{\w_i, \w_j}^{2}\w_j	-  \inner{\w_i, 
		\e_{d-1}}^{2}\e_{d-1}}.
\end{align*}
The term in the above corresponding to $i=d-1$ is 
$\prn*{(1-t^3)^{2/3}}^2(1-t^3)^{1/3} + 
\prn*{t(1-t^3)^{1/3}}^2t-(1-t^3)^{2/3}=0$. Similarly, for $i=d$, 
$t^5 + \prn*{(1-t^3)^{1/3}t}^2(1-t^3)^{1/3} - t^2 = 0$. Moreover, by 
\pref{eqn:obj},
\begin{align*}
\ploss(\fC_{5,t}(d)) 
&= \sum_{i,j=1}^{d} \inner{\w_i, \w_j}^3
-2\sum_{i=1}^{d}  \sum_{j=1}^{d} \inner{\w_i, \e_j}^3
+ \sum_{i,j=1}^{d} \inner{\e_i, \e_j}^3\\
%&= 
%d-2 + 2(t(1-t^3)^{1/3})^3 + t^6 + (1-t^3)^{2}
%-2(d-2 + t^3 + 1-t^3)
%+ d\\
&= 2t^3(1-t^3) + t^6 + 1-2t^3 + t^6 = 1.
\end{align*}
Therefore, the level set $\ploss^{-1}(1)$ is unbounded, hence $\ploss$ 
is not proper.

\subsection{\pfref{prop:inv}}

Let $\sigma = (\sigma_1,\sigma_2)\in 
(S_k)_{\balpha}\times \isov_{\bbeta,V} =
(S_k)_{\balpha}\times \Pi_2((S_h\times 
S_d)_{V}\cap 
((S_h)_{\bbeta}\times S_d) )$. 
Since $\sigma_2\in \isov_{\bbeta,V}$, there 
exists $\rho\in (S_h)_{\bbeta}$ 
such that $(\rho,\sigma_2)\in (S_h\times 
S_d)_{V}$. In particular, 
$\rho^{-1}\in(S_h)_{\bbeta}$ and 
$(\rho^{-1},\sigma_2^{-1})\in 
(S_h\times 
S_d)_{V}$. By \pref{prob:kernel} and the 
permutation 
invariance
of $\ker$,
\begin{align*} \label{kernel_inv}
\ploss_\ker&(\sigma W;V) \\
&=\sum_{i,j=1}^{k} 
\alpha_i\alpha_j\ker\prn{\sigma_2\w_{\sigma_1^{-1}(i)}, 
\sigma_2\w_{\sigma_1^{-1}(j)}}
-2\sum_{i=1}^{k} \sum_{j=1}^{h} 
\alpha_i\beta_j\ker\prn{\sigma_2\w_{\sigma_1^{-1}(i)}, 
\v_j}
+ \sum_{i,j=1}^{h} \beta_i\beta_j\ker\prn{\v_i, \v_j}\\
%&=
%\sum_{i,j=1}^{k} 
%\alpha_{i}\alpha_j\ker\prn{\w_{\sigma_1^{-1}(i)}, 
%\w_{\sigma_1^{-1}(j)}}
%-2\sum_{i=1}^{k} \sum_{j=1}^{h} 
%\alpha_i\beta_j\ker\prn{\w_{\sigma_1^{-1}(i)}, 
%	\sigma_2^{-1}\v_j}+ \sum_{i,j=1}^{h} 
%	\beta_i\beta_j\ker\prn{\v_i, 
%	\v_j}\\
&=
\sum_{i,j=1}^{k}\alpha_{\sigma_1(i)}\alpha_{\sigma_1(j)} 
\ker\prn{\w_{i}, \w_{j}}
-2\sum_{i=1}^{k} \sum_{j=1}^{h} 
\alpha_{\sigma_1(i)}\beta_j\ker\prn{\w_{i}, 
\sigma_2^{-1}\v_j}+ \sum_{i,j=1}^{h} \beta_i\beta_j\ker\prn{\v_i, 
\v_j}\\
%&=
%\sum_{i,j=1}^{k}\alpha_{i}\alpha_{j} 
%\ker\prn{\w_{i}, \w_{j}}
%-2\sum_{i=1}^{k} \sum_{j=1}^{h} 
%\alpha_{i}\beta_j\ker\prn{\w_{i}, 
%	\sigma_2^{-1}\v_j}+ \sum_{i,j=1}^{h} 
%	\beta_i\beta_j\ker\prn{\v_i, 
%	\v_j}\\
\end{align*}
\begin{align*}
&=
\sum_{i,j=1}^{k}\alpha_{i}\alpha_{j} 
\ker\prn{\w_{i}, \w_{j}}
-2\sum_{i=1}^{k} \sum_{j=1}^{h} 
\alpha_{i}\beta_{\rho^{-1}(j)}\ker\prn{\w_{i}, 
\sigma_2^{-1}\v_j} + \sum_{i,j=1}^{h} 
\beta_{\rho^{-1}(i)}\beta_{\rho^{-1}(j)}\ker\prn{\v_i,
\v_j}\\
%&=
%\sum_{i,j=1}^{k}\alpha_i\alpha_j \ker\prn{\w_{i}, 
%\w_{j}}
%-2\sum_{i=1}^{k} \sum_{j=1}^{h} 
%\alpha_i\beta_j\ker\prn{\w_{i}, 
%	\sigma_2^{-1}\v_{\rho(j)}}+ \sum_{i,j=1}^{h} 
%\beta_i\beta_j
%	\ker\prn{\v_{\rho(i)}, \v_{\rho(j)}}\\
&=
\sum_{i,j=1}^{k} \alpha_i\alpha_j\ker\prn{\w_{i}, \w_{j}}
-2\sum_{i=1}^{k} \sum_{j=1}^{h} \alpha_i\beta_j\ker\prn{\w_{i}, 
\sigma_2^{-1}\v_{\rho(j)}}+ \sum_{i,j=1}^{h} \beta_i\beta_j
\ker\prn{\sigma_2^{-1}\v_{\rho(i)}, \sigma_2^{-1}\v_{\rho(j)}}\\	
&= \ploss(W;(\rho^{-1},\sigma_2^{-1})V)\\
&= \ploss(W;V),
\end{align*}
concluding the proof.

\section{Existence and uniqueness of Puiseux series} 
\label{sec:puiseux}
We describe a construction of Puiseux series using the (real or complex) analytic implicit function theorem \cite{ArjevaniField2022annihilation,ArjevaniField2020}, see generally\cite{wall2004singular,de2013local,mcdonald2002fractional} for additional approaches. The method is illustrated with reference to $\fc_{\bC}$ (see \pref{rem:complex}) with isotropy group $\Delta S_{d}$ and so given by two Puiseux series $\xi_1(d)$ and $\xi_2(d)$. The system of critical equations restricted to the fixed point space $M(d,d)^{\Delta S_d}$ is given in (\ref{deltasd_cubic}); denote it by $E(\xi_1,\xi_2,d)$. Repeatedly applying a Newton polyhedron argument (see \pref{sec:rank_for}) gives
\begin{align*}
\xi_1(d) &= 1 - \frac{2 i}{3 
d^{\frac{2}{4}}} + 
\frac{\sqrt{2} \cdot \left(1 + 
i\right)}{d^{\frac{3}{4}}} - 
\frac{77}{36 d} + \frac{43 
\sqrt{2} \cdot \left(1 - i\right)}{96 
d^{\frac{5}{4}}} + \frac{2821 i}{5184 
d^{\frac{6}{4}}}   \\&+ \frac{8491 \sqrt{2} \left(-1 
- 
i\right)}{18432 d^{\frac{7}{4}}} + O(d^{-2}),\\
\xi_2(d) &= \frac{\sqrt{2} \cdot \left(1 + 
i\right)}{2 
d^{\frac{3}{4}}} - \frac{5}{4 d}
+ \frac{49 
\sqrt{2} \cdot \left(1 - 
i\right)}{192 d^{\frac{5}{4}}}
+ \frac{119 i}{96 d^{\frac{6}{4}}} 
+ \frac{39799 \sqrt{2} \left(-1 - i\right)}{36864 
d^{\frac{7}{4}}}\\& + \frac{13}{9 
d^{2}}  + O(d^{-\frac{9}{4}}).
\end{align*}
The trimmed Puiseux 
series are defined by 
\begin{align*}
\xi_1^{(n_1)}(d) &\coloneqq \sum_{i=0}^{n_1} 
[d^{\frac{-i}{4}}]\xi_1(d),\\
\xi_2^{(n_2)}(d) &\coloneqq \sum_{i=0}^{n_2} 
[d^{\frac{-i}{4}}]\xi_2(d),
\end{align*}
where $[d^\alpha]P$ denotes the $\alpha$-order 
term of a Puiseux series $P$. Thus, for example, 
$\xi_1^{(2)}(d) 
= 1 - \frac{2 
i}{3 d^{\frac{2}{4}}}$ and $\xi_2^{(2)}(d) = 0$.
Next, set $s = d^{\frac{-1}{4}}$. This 
substitution 
turns $\xi_i(d)$ into a power series having {integral }exponents. The set of equations so obtained (divided by 2) is
\begin{align*}
\overline{E}(x,y,s) \defeq E(\xi_1^{(n_1)}(s) + 
s^{n_1+1}x, \xi_2^{(n_2)}(s) + s^{n_2+1}y, s^{-4}),
\end{align*}
where the trimming parameters, $n_1$ and $n_2$, 
are chosen so that the Jacobian 
of $\overline{E}$,
\begin{align*}
\begin{pmatrix}
\frac{\partial \overline{E}_1}{\partial x}, 
\frac{\partial \overline{E}_1}{\partial y}\\
\frac{\partial \overline{E}_2}{\partial x},
\frac{\partial \overline{E}_2}{\partial y}
\end{pmatrix},
\end{align*}
at the {initial} conditions
$x = [s^{n_1+1}]\xi_1(s)$, $y= 
[s^{n_2+1}]\xi_2(s)$ and $s=0$
is non-singular (assuming momentarily that powers of $s$ occurring in each equation have been canceled). We may now invoke the analytic 
implicit function theorem (here, in the complex 
category) and so present $x(s)$ and 
$y(s)$ as power series in $s$ in a 
neighborhood of $s=0$. Reversing the substitution, 
that is, substituting $d^{\frac{-1}{4}}$ for $s$, 
completes the derivation (note, convergence is now in 
a neighborhood of $d=\infty$, that is, for 
sufficiently 
large $d$). For $\fc_{\bC}$, suffices it to take 
$n_1= 1$ and 
$n_2 = 2$. Thus, $\overline{E}(x,y,s) = 
E(1 + s^2 x, s^3 y, s^{-4})$. We have, for 
example,
\begin{align*}
\overline{E}(x,y,s)_1 &= - 12 
s^{15} y^{5} + 27 s^{14} x y^{4} - 12 s^{13} x^{2} 
y^{3} + s^{12} \left(- 6 x^{3} y^{2} + 27 
y^{4}\right) \\&+ s^{11} \left(- 24 x y^{3} + 24 
y^{5}\right) + s^{10} \left(3 x^{5} - 18 x^{2} 
y^{2} - 42 x y^{4}\right) \\&+ s^{9} 
\left(12 
x^{2} y^{3} - 12 y^{3}\right) + s^{8} 
\left(15 
x^{4} + 6 x^{3} y^{2} - 18 x y^{2} - 42 y^{4}\right) 
\\&+ s^{7} \left(24 x y^{3} - 15 
y^{5}\right) + s^{6} \left(30 x^{3} + 18 x^{2} y^{2} + 15 x 
y^{4} - 6 y^{2}\right) \\&+ 12 s^{5} y^{3} + 
s^{4} \left(27 x^{2} + 18 x y^{2} + 15 y^{4}\right) + 
3 s^{3} y^{5} + s^{2} \left(9 x + 6 
y^{2}\right).
\end{align*}
It follows that $s^2$ should be canceled out from $\overline{E}(x,y,s)_1$. Similarly, one 
deduces that $s^3$ should be canceled out from 
$\overline{E}(x,y,s)_2$. 
Let 
$\tilde{E}(x,y,s)$ denote the resulting set of 
equations. 
Plug-in  $x 
=[s^{n_1+1}]\xi_1(s) = \frac{-2i}{3} ,~y= 
[s^{n_2+1}]\xi_2(s) = \frac{\sqrt{2}}{2} + 
\frac{\sqrt{2} i}{2}$ and $s=0$ into $\tilde{E}$ 
and observe that the initial conditions hold and 
that the Jacobian with respect to $x,y$ 
(evaluated at the initial conditions) 
\begin{align*}
J = \left[\begin{matrix}9 & 6 \sqrt{2} + 6 
\sqrt{2} i\\0 & 15 \left(\frac{\sqrt{2}}{2} + 
\frac{\sqrt{2} i}{2}\right)^{4} + 
3\end{matrix}\right]	
\end{align*}
is indeed non-singular, concluding the derivation. 
Although only the first few terms are used, by uniqueness, the rest of the terms obtained by the Newton polyhedron argument are correct.

\end{document}